\documentclass[12pt]{article}
\usepackage{bbm}
\usepackage{mathrsfs}
\usepackage{amsfonts}
\usepackage{amsmath,amssymb,latexsym,amsfonts}

\usepackage{CJK,color}
\usepackage{amssymb}
\usepackage{amsmath}
\usepackage{graphics}
\usepackage{epsfig}
\usepackage{indentfirst}
\usepackage{cite}

\definecolor{red}{rgb}{1.00,0.00,0.00}
\definecolor{blue}{rgb}{0.00,0.00,0.63}
\definecolor{black}{rgb}{0.00,0.00,0.00}

\newtheorem{claim}{\bf \t}[part]

\newtheorem{Lemma}{Lemma}[part]
\newtheorem{Proposition}{Proposition}[part]
\newtheorem{Remark}{Remark}[part]

\newtheorem{Theorem}{Theorem}[part]

\numberwithin{Assumption}{section} \numberwithin{Corollary}{section}
\numberwithin{Definition}{section} \numberwithin{equation}{section}
\numberwithin{Example}{section} \numberwithin{Lemma}{section}
\numberwithin{Proposition}{section} \numberwithin{Remark}{section}
\numberwithin{Theorem}{section}

\def\t{\theta}

\def \qed{{\hbox{ }\hfill$\Box$}}

\def\text#1{{\rm #1}}

\begin{document}
\date{}

\title{\Large \bf  Formation of Singularities of Spherically Symmetric Solutions  to the 3D Compressible Euler Equations and Euler-Poisson Equations }

\author{Hai-Liang Li\  and\ Yuexun Wang}

\maketitle
\begin{abstract}
By introducing a new averaged quantity with a fast decay weight to perform Sideris's argument \cite{sideris1985formation} developed for the Euler Equations,  we extend the formation of singularities of classical solution to the 3D Euler Equations established in \cite{sideris1985formation,makino1986solution} for the initial data with compactly supported disturbances to the spherically symmetric solution with general initial data in Sobolev space. Moreover, we also prove the formation of singularities of the spherically symmetric solutions to the 3D Euler-Poisson Equations, but remove the compact support assumptions on the initial data in \cite{makino1990solution,perthame1990non}. Our proof also simplifies that of \cite{Lei2013singularities} for the Euler Equations and is undifferentiated in dimensions.  
\end{abstract}

\setcounter{equation}{0}
\setcounter{Assumption}{0}
\setcounter{Theorem}{0}
\setcounter{Proposition}{0}
\setcounter{Corollary}{0}
\setcounter{Lemma}{0}
\mbox{}

\section{Introduction}
It is well-known that Sideris \cite{sideris1985formation} first proved the important fact that the \(C^1\) solution to the three-dimensional Euler equations for compressible fluids must develop singularity in finite time. His proof consists of two  critical ingredients, one is the the finite propagation  of compactly supported disturbances of  the solution, the other is the evolution of certain averaged quantities formed out of the solution. In the proof, he assumes that the initial density \(\rho_0\) is  positive on the entire space (\(\inf \rho_0>0\)) and has a positive background state outside  a ball (\(\rho_0=\bar{\rho}, B_{R}^c\)), and the initial velocity \(u_0\) is compactly supported on the same ball. This design makes sure that, on the one hand, the Euler equations may be written as a
 positive definite, symmetric hyperbolic system by the usual symmetrization, and thus possess a unique, local \(C^1\) solution \((\rho-\bar{\rho}, u)\) via the theory \cite{kato1975cauchy,majda2012compressible}, on the other hand, by the local energy estimates \cite{sideris1984formation}, the solution \((\rho-\bar{\rho}, u)\) is supported  on a ball, whose radius is determined by the sound speed corresponding to the background state \(\bar{\rho}\), and thus has linear spreading rate. To prove the first main result \cite{sideris1985formation}, Sideris introduced an averaged quantity, the radial momentum \(\int \rho ux\,\mathrm{d} x\), whose  evolution  connects with some other averaged quantities, for example, the mass,  energy and moment of inertia \(\int \rho |x|^2\, \mathrm{d} x\) (refer to \cite{sideris2014spreading} for an exhaustive account). The crucial idea is to estimate the upper bound of the moment of inertia which  heavily relies on the finite propagation speed. The sign of relative total pressure mainly depends on the  convexity of the pressure but may be independent of the finite propagation speed (see \cite{sideris2003long}).
Finally, the radial momentum obeys a Riccati type inequality which forces the lifespan of the solution  to be finite if the initial velocity is supersonic in some region. By studying other alternative averaged quantities, in the third main result,  Sideris established the formation of singularities without any condition of largeness on the initial velocity, which also provides an upper bound of the lifespan in exponential type \(\exp(C/\varepsilon^2)\), if the initial data has a perturbation of order \(\varepsilon\) from the positive background state. In another paper \cite{sideris1991lifespan},  Sideris  refined the upper bound \(\exp(C/\varepsilon^2)\) of the lifespan to the order \(\exp(C/\varepsilon)\) for irrotational initial velocity, which was based on  a key observation  that the linear operator of acoustics is invariant under the Lorentz transformations for irrotational velocity fields.

As the initial density is surrounded by vacuum, it does not seem yet to be straightforward to apply the theory \cite{kato1975cauchy,majda2012compressible}  directly  to the Euler equations via  the usual symmetrization for lack of positive definiteness.  Makino, Ukai and Kawashima \cite{makino1986solution} made use of the sound speed  as an independent variable in place of the density which together with the original velocity  symmetrized the Euler equations to a positive definite hyperbolic system, and thus the well-posedness of classical solution follows. Then, under the further compact support assumptions on the initial density and velocity, the authors showed that the support of the solution does not change  with respect to time as long as the solution is tame. This fact gives rise to the the moment of inertia an upper bound independent of the time, however, whose lower bound spreads quadratically with respect to time. This obvious contradiction implies that the blowup must occur in finite time for the tame solution. Makino and Ukai \cite{makino1987existence} continually  applied the same symmetrization \cite{makino1986solution} to establish  the local solvability of three-dimensional Euler-Poisson equations with the gravitational forces by a fixed point argument, the idea in which is also effective in the presence of the electrostatic forces.

If the initial data are spherically symmetric and have compact support, for the three-dimensional Euler-Poisson equations with gravitational forces, Makino and Perthame \cite{makino1990solution} proved that the support of the tame solution does not spread over time with the aid of the sign of the gravitational force, which presents a similar blowup result of classical to \cite{makino1986solution}.
When the external forces are the electrostatic force, Perthame \cite{perthame1990non} showed that the support of spherically symmetric tame solution grows quadratically. This means that the moment of inertia is sandwiched between a pair of quadratic functions whose coefficients consist of the total energy, mass and so on. By comparing the coefficients of the two quadratic functions, the blowup of classical solution follows if the initial energy is large in comparison with the initial mass.

In \cite{sideris1985formation}, the maximum propagation speed of a front into a constant state is a-priorily determined by the positive background state. In \cite{makino1986solution,makino1990solution,perthame1990non}, the maximum propagation speed of compact support is also a-priorily determined by the initial data.  However it is not clear how the support of the general spherically symmetric solution in Sobolev space propagates, especially, as the presence of non-local terms in the compressible Euler-Poisson equations.  To overcome the difficulty caused by the lack of the finite speed propagation, one needs to consider some new averaged quantities which may localize the far field for applying  Sideris's argument \cite{sideris1985formation}. Starting  from the formulation that,  as observed in \cite{sideris1985formation}, the Euler equations may be written as a quasi-linear wave-type equation in terms of density with inhomogeneous terms involving the velocity, Lei, Du and Zhang \cite{Lei2013singularities} studied the evolutions of some new averaged quantities with fast decay weights which led to some formation of singularities of spherically symmetric
 solutions to the Euler equations in two and three spatial dimensions, in which the authors
 removed the compactly supported disturbances assumptions on the initial data. The proofs present delicate differences between the two and three spatial dimensions due to the different weights.

There are extensive studies on the formation of singularities of the one-dimensional Euler equations and Euler-Poisson equations which are most closely related to the method of characteristics, also developed at the background of hyperbolic systems of
conservation laws in one spatial dimension. However we cannot exhaust all of them and here only list some of them for the reader's reference convenience \cite{Dafermos1985development, Lax1964development,Liu1979development, Wang1998formation}, one may also refer to \cite{Chen2002Cauchy,Dafermos2010Hyperbolic} for a fairly complete list of the references.

This work in the present paper  is  devoted to the investigation of the formation of singularities of spherically symmetric solutions with general initial data in Sobolev space to the three-dimensional Euler equations and Euler-Poisson equations with gravitational force or electrostatic force.
However it seems that the method of \cite{Lei2013singularities} doesn't apply directly to the Euler-Poisson equations. Different from \cite{Lei2013singularities}, we study the equation of the velocity and work with spherical coordinates instead of Eulerian coordinates, which leads us to find a very simple weight \(\exp(-r)\) that is effective to both the Euler equations and Euler-Poisson equations. On the other hand,  our proof simplifies that of \cite{Lei2013singularities} and is also undifferentiated between the two-and three-dimensional Euler equations. Specifically, we introduce a simple averaged quantity \(F(t)\) with the new weight \(\exp(-r)\),  which is efficient in performing Sideris's argument \cite{sideris1985formation} to the Euler-Poisson equations, as well as the Euler equations.  Due to the good sign of the gravitational forces,  we may handle the Euler equations and Euler-Poisson equations with gravitational forces at the same time. Following Sideris's argument \cite{sideris1985formation}, we study the evolution of \(F(t)\) that yields a middle term, the weighted kinetic energy \(Q(t)\), which may bound \(F(t)\) from above, but independent of the time.   On the occasions of the Euler-Poisson equations with electrostatic forces, no good sign to use, we have to employ the weighted Hardy inequality to treat with the high singularity near \(r=0\). However one still can control the other middle term coming from  the electrostatic forces, but in terms of the initial mass and physical energy instead, and all the estimates are independent of the time. Finally \(F(t)\) satisfies a Riccati type inequality and thus the singularity forms in finite time. The above formation of singularities of the Euler equations and Euler-Poisson equations with gravitational forces provide an upper bound for the lifespan to a time of the order \(C/\varepsilon\), this together with the lower bound \(C/\varepsilon\) from \cite{makino1986solution,makino1987existence}, yields a by-product, the sharp bound of the lifespan in these two situations.

 Section 2 contains the  reformulation of
 the problem, as well as some necessary preliminaries. In Section 3 we revisit the local well-posedness of classical solutions to the three-dimensional Euler equations and Euler-Poisson equations with gravitational forces or electrostatic forces. We establish the formation of singularities for spherically symmetric solutions in Sobolev space in Section 4. Section 5 presents some examples to justify the conditions in the formation of singularities.

\section{Reformulation of the problem}
We consider here the  Euler equations for
an isentropic ideal fluid in  three-dimension which are governed by

\begin{eqnarray}\label{1}
	\left\{ \begin{array}{ll}
		\partial_t\rho+\textrm{div}(\rho
		u)=0,\\
		\partial_t(\rho u)+\textrm{div}(\rho u\otimes u)+\nabla p=0.
	\end{array}
	\right.
\end{eqnarray}
Here  \(\rho,u=(u_1,u_2,u_3)\) and \(p=p(\rho)\) denote the
density, velocity and pressure, respectively.  For the ideal polytropic gas, the equation of state is given by \(p=A\rho^\gamma\) with the adiabatic index \(\gamma> 1\) and a constant \(A>0\).
When considering gas motion under self-gravitation or describing the dynamics of a plasma,  these physical
phenomenon may be modeled by the Euler-Poisson system
\begin{eqnarray}\label{2}
	\left\{ \begin{array}{ll}
		\partial_t\rho+\textrm{div}(\rho
		u)=0,\\
		\partial_t(\rho u)+\textrm{div}(\rho u\otimes u)+\nabla p=\pm\rho \nabla \phi,\\
		\Delta \phi=4\pi \mathcal{G}\rho,
	\end{array}
	\right.
\end{eqnarray}
where we use the  signs \("\pm"\) to distinguish  electrostatic forces \("+"\) and self-gravitational forces \("-"\). We set both the constant \(A\) and universal gravitational constant \(\mathcal{G}\) to be unit for simplicity.
One can solve  the Newtonian potential \(\phi\) from the elliptic equation in \eqref{2} in terms of  \(\rho\) to reduce \eqref{2} to a much simpler system
\begin{eqnarray}\label{3}
	\left\{ \begin{array}{ll}
		\partial_t\rho+\textrm{div}(\rho
		u)=0,\\
		\partial_t(\rho u)+\textrm{div}(\rho u\otimes u)+\nabla p=\pm \rho \nabla G\ast\rho,\\
	\end{array}
	\right.
\end{eqnarray}
in which  \(G\) stands for the three-dimensional Green's function (factoring out a constant) which is given by \(G(x)=-1/|x|\). We denote by \(s(\rho)\) the local sound speed which is determined by the formula \(s(\rho)=\sqrt{p_\rho(\rho)}\). Denote \( \zeta=\frac{2}{\gamma-1}s(\rho)\),  as an independent variable instead of \(\rho\) together with velocity, which may symmetrize the Euler equations \eqref{1} and Euler-Poisson equations \eqref{3} to
\begin{align}\label{4}
	\partial_tV+\sum_{i=1}^nA_{i}(V)\partial_{x_i}V+\delta B(V)=0.
\end{align}
Here \(\delta\) takes the values \(0,1,-1\) which are corresponding to the Euler equations, Euler-Poisson equations with electrostatic forces and with gravitational forces, respectively. The notations in the system \eqref{4} are specifically as follows
\[
V =
\begin{pmatrix}
\zeta \\
u
\end{pmatrix}, \
B(V) =
\begin{pmatrix}
0 \\
-\big[\frac{(\gamma-1)^2}{4\gamma}\big]^\frac{1}{\gamma-1}\nabla G\ast\zeta^{\frac{2}{\gamma-1}}
\end{pmatrix}
,\]
\[
A_1(V)=
\begin{pmatrix}
u_1 & \frac{(\gamma-1)\zeta}{2}&0&0 \\
\frac{(\gamma-1)\zeta}{2} & u_1&0&0\\
0&0&u_1&0\\
0&0&0&u_1
\end{pmatrix},
\]

\[
A_2(V)=
\begin{pmatrix}
u_2 & 0&\frac{(\gamma-1)\zeta}{2}&0 \\
0 & u_2&0&0\\
\frac{(\gamma-1)\zeta}{2}&0&u_2&0\\
0 & 0&0&u_2\\
\end{pmatrix}
,\]

\[
\ A_3(V)=
\begin{pmatrix}
u_3 & 0&0&\frac{(\gamma-1)\zeta}{2} \\
0 & u_3&0&0\\
0&0&u_3&0\\
\frac{(\gamma-1)\zeta}{2} & 0&0&u_3\\
\end{pmatrix}.
\]

For notational convenience, we define  the partial and total energy functionals as
\begin{align*}
	\mathcal{E}^{(k)}(t) =\|V^{(k)}(t,\cdot)\|_{L^2(\mathbb{R}^3)}^2
\end{align*}
and
\begin{align*}
	\mathcal{E}(t) =\sum_{k=0}^3\mathcal{E}^{(k)}(t),
\end{align*}
respectively. Here the maximum index in the above summation is to make sure that the solutions constructed under these energy functionals are classical via Sobolev embedding. We denote by \(\mathcal{X}_T\) the Sobolev space
\[C\big([0,T);H^3(\mathbb{R}^3)\big)\cap C^1\big([0,T);H^2(\mathbb{R}^3)\big)\]
for simplicity which will be frequnetly used in the well-posedness and singularity formation theories.

We are interested in studying formation of singularities of classical solutions with conservative mass and physical energy to the three-dimensional Euler equations and  Euler-Poisson equations with gravitational forces  or electrostatic forces. We use \(M(t)\) and \(E_\delta(t)\) to denote the mass  and  physical energy respectively which are defined by
\begin{align*}
	M(t)=\int_{\mathbb{R}^3}\rho\,\mathrm{d} x
\end{align*}
and
\begin{align*}
	E_\delta(t)=\int_{\mathbb{R}^3}\big(\frac{1}{2}\rho u^2+\frac{1}{\gamma-1}\rho^\gamma-\frac{\delta}{2}\rho G\ast\rho\big)\,\mathrm{d} x.
\end{align*}
Here \(\delta\) has the same meaning as \eqref{4}. We should mention that \(E_1(t)\)  is always nonnegative  due to  the negative sign of \(G(x)\). However  \(E_{-1}(t)\) may be nonnegative or negative, whose sign depends on initial setup.

The notation \(C\)  always denotes a nonnegative universal constant which may be different from line to line but is
independent of the parameters involved. Otherwise, we will specify it by  the notation \(C(a,b,\dots)\).
We write \(f\lesssim g\) (\(f\gtrsim g\)) when \(f\leq  Cg\) (\(f\geq  Cg\)), and \(f \eqsim g\) when \(f \lesssim g \lesssim f\). The notation \(\varepsilon\) always stands for a  sufficiently small positive number throughout the paper.

\section{Local Well-Posedness}
In this section we will revisit the well-posedness of classical solutions to the Euler equations and  Euler-Poisson equations with electrostatic forces or gravitational forces established in \cite{makino1986solution,makino1987existence}. To avoid the trivial cases we always assume in the present paper that the initial data satisfy
\[
\rho_0\geq0,\ 0<\mathcal{E}(0),M(0),E_\delta(0)<\infty.
\]

We first state the well-posedness of the Euler equations (\(\delta=0\)).
\begin{Proposition}\label{prop:Euler}  If \(1<\gamma\leq\frac{5}{3}\).  Then there exist a positive number
	\begin{align}\label{5}
	T\gtrsim\mathcal{E}^{-\frac{1}{2}}(0)
	\end{align}
	and a unique  solution \((\rho,u)\in \mathcal{X}_T\) of the Euler equations \eqref{1}  with  \((\rho,u)|_{t=0}=(\rho_0,u_0)\) satisfying \(\rho\geq0\), which possesses the conservations of total mass and energy.
\end{Proposition}

\textbf{Proof.} Since the  system \eqref{4} is a positive definite hyperbolic system,  the standard theory \cite{kato1975cauchy,majda2012compressible} applies to show that the system \eqref{4} is  locally well-posed in \(\mathcal{X}_T\) for some positive time \(T\).
Moreover  \(\mathcal{E}(t)\) satisfies the energy estimate
\begin{align*}
\frac{\mathrm{d}}{\mathrm{d} t}\mathcal{E}(t)\lesssim \mathcal{E}^{\frac{3}{2}}(t),
\end{align*}
which  provides a lower bound \(T\gtrsim\mathcal{E}^{-\frac{1}{2}}(0)\). One solves the first equation of \eqref{4} to deduce \(\zeta\geq0\).

Since \(\zeta\) belongs to \(\mathcal{X}_T\), so does \(\rho\) if \(1<\gamma\leq\frac{5}{3}\). This is a consequence of \(\rho=C(\gamma)\zeta^\frac{2}{\gamma-1}\) (the definition of the sound speed).  Then  Sobolev embedding 	
yields \((\rho,u)\in C^1([0,T)\times \mathbb{R}^3)\), which implies that the inverse map of \((\rho,u)\mapsto(\zeta,u)\) is continuously differentiable. So \((\rho,u)\) solves the Cauchy problem of \eqref{1}.
Since we have shown \((\rho,u)\in \mathcal{X}_T\), which  guarantees the integration by parts  in  proving the conservations of mass.
\qed

The well-posedness of the Euler-Poisson equations (\(\delta=\pm 1\)) may be stated as follows.
\begin{Proposition}\label{prop:EPEG} If \(1<\gamma\leq\frac{5}{3}\). Then there exist a positive number
	\begin{align}\label{5.5}
	T\gtrsim
	\max\{\mathcal{E}^{-\frac{1}{2}}(0),\mathcal{E}^{\frac{\gamma-3}{2(\gamma-1)}}(0)\}
	\end{align}
	and a unique  solution \((\rho,u)\in \mathcal{X}_T\) of the Euler-Poisson equations \eqref{3}  with   \((\rho,u)|_{t=0}=(\zeta_0,u_0)\) satisfying \(\rho\geq0\), which possesses the  conservations of total mass  and energy.
\end{Proposition}
\textbf{Proof.} We here give only an a priori estimate. Applying \(D^k\), \(0\leq k\leq3\) to \eqref{4}, then multiplying it by \(D^kV\) and using integration by parts to get
\begin{align}\label{6}
\frac{\mathrm{d}}{\mathrm{d} t}\mathcal{E}(t)\lesssim \mathcal{E}^{\frac{3}{2}}(t)+\sum_{k=0}^3\|D^k(\nabla G\ast\zeta^{\frac{2}{\gamma-1}})\|_{L^2(\mathbb{R})}\|D^kU\|_{L^2(\mathbb{R})}.
\end{align}

Note that \(|\nabla G|=1/|x|^2\), we employ  Hardy-Littlewood-Sobolev inequality from Appendix to deduce  	
\begin{equation}\label{7}
\begin{aligned}
\|D^k(\nabla G\ast\zeta^{\frac{2}{\gamma-1}})\|_{L^2(\mathbb{R}^3)}
&= \|\nabla G\ast D^k(\zeta^{\frac{2}{\gamma-1}})\|_{L^2(\mathbb{R}^3)}
\lesssim\|D^k(\zeta^{\frac{2}{\gamma-1}})\|_{L^{\frac{6}{5}}(\mathbb{R}^3)}\\
&\lesssim\|D^k(\zeta^{\frac{2}{\gamma-1}})\|_{L^1(\mathbb{R}^3)}^{\frac{2}{3}}\|D^k(\zeta^{\frac{2}{\gamma-1}})\|_{L^2(\mathbb{R}^3)}^{\frac{1}{3}}.
\end{aligned}
\end{equation}
The Gagliardo-Nirenberg  inequality shows 	
\begin{align}\label{8}
\|D^k(\zeta^{\frac{2}{\gamma-1}})\|_{L^1(\mathbb{R}^3)}\lesssim\|\zeta^{\frac{2}{\gamma-1}}\|_{L^1(\mathbb{R}^3)}^{\frac{3}{3+2k}}\|D^k(\zeta^{\frac{2}{\gamma-1}})\|_{L^2(\mathbb{R}^3)}^{\frac{2k}{3+2k}}.
\end{align}
Then substituting \eqref{8} to \eqref{7} yields
\begin{align}\label{9}
\|D^k(\nabla G\ast\zeta^{\frac{2}{\gamma-1}})\|_{L^2(\mathbb{R}^3)}
\lesssim\|\zeta^{\frac{2}{\gamma-1}}\|_{L^1(\mathbb{R}^3)}^{\frac{2}{3+2k}}\|D^k(\zeta^{\frac{2}{\gamma-1}})\|_{L^2(\mathbb{R}^3)}^{\frac{1+2k}{3+2k}}.
\end{align}

Note that  \(1<\gamma\leq\frac{5}{3}\),  each term \(D^k(\zeta^{\frac{2}{\gamma-1}})\), \(k=0,1,2,3\) may be bounded by the norm \(H^3(\mathbb{R}^3)\) of \(\zeta\) via Sobolev embedding. Thus one can estimate \eqref{9} continually to obtain
\begin{align}\label{10}
\|D^k(\nabla G\ast\zeta^{\frac{2}{\gamma-1}})\|_{L^2(\mathbb{R}^3)}
\lesssim\|\zeta^{\frac{4-2\gamma}{\gamma-1}}\|_{L^\infty(\mathbb{R}^3)}^{\frac{2}{3+2k}}\|\zeta^2\|_{L^1(\mathbb{R}^3)}^{\frac{2}{3+2k}}
\|\zeta\|_{H^3(\mathbb{R}^3)}^{\frac{2(1+2k)}{(\gamma-1)(3+2k)}}
\lesssim\|\zeta\|_{H^3(\mathbb{R}^3)}^{\frac{2}{\gamma-1}}.
\end{align}
Inserting \eqref{10} into \eqref{6} gives
\begin{eqnarray*}
	\frac{\mathrm{d}}{\mathrm{d} t}\mathcal{E}(t)\lesssim
	\mathcal{E}^{\frac{3}{2}}(t)+\mathcal{E}^{\frac{\gamma+1}{2(\gamma-1)}}(t),
\end{eqnarray*}
which in view of Gr\"{o}nwall's inequality yields \eqref{5.5}. The rest is akin to the corresponding part in the proof of Proposition \ref{prop:Euler}.

\qed

\section{Formation of Singularities}

In this section we study the formation singularities of  spherically symmetric solutions to the three-dimensional Euler equations and  Euler-Poisson equations with gravitational forces  or electrostatic forces. For spherically symmetric motions, the solution to \eqref{1} or \eqref{3} established in Section 3, emanating from spherically symmetric initial data
\begin{align}\label{10.5}
\rho_0(x)=\rho_0(r),\ u_0(x)=\frac{x}{r}v_0(r),\ r=|x|,
\end{align}
 possess the spherically symmetric form
\begin{align}\label{11}
\rho(t,x)=\rho(t,r),\ u(t,x)=\frac{x}{r}v(t,r),
\end{align}
which also satisfies
\begin{eqnarray}\label{12}
\left\{ \begin{array}{ll}
\partial_t\rho+\frac{1}{r^2}\partial_r(r^2\rho
v)=0,\\
\rho(\partial_tv+v\partial_rv)+\partial_rp=\delta \frac{4\pi\rho }{r^2}\int_0^r\rho(t,l) l^2\,\mathrm{d}l.
\end{array}
\right.
\end{eqnarray}
Here \(\delta\) has the same meaning with \eqref{4} and \(r>0\). Let \((\rho,u)\in \mathcal{X}_T\). As observed in \cite{Lei2013singularities}, we see from \eqref{11} that  each component of \(u\) is an odd function, so the velocity \(u\) should satisfy \(u(t,0)=0\) which in turn implies \(v(t,0)=0\). Since the density \(\rho\) is radial, thus the characteristic curve always maps the center \(r=0\) to itself as long as the solution is smooth, or there must be at least two characteristic curves intersecting at \(r=0\). It follows from the mass equation of \eqref{1} or \eqref{3} that \(\rho(t,0)=0\) if the initial density is imposed on \(\rho_0(0)=0\). In the following we will study formation singularities of  spherically symmetric solutions with the initial data \(\rho_0(0)=0\). For convenience, we formulate the above observation as
\begin{align}\label{12.5}
\rho(t,0)=0,\ v(t,0)=0.
\end{align}

The mass and physical energy defined in Section 3 may be transfered into the following corresponding to spherically symmetric solutions to \eqref{12} as
\begin{align*}
M(t)=4\pi\int_0^\infty\rho(t,r) r^2\,\mathrm{d} r
\end{align*}
and
\begin{align*}
E_\delta(t)&=4\pi\int_0^\infty\big(\frac{1}{2}\rho v^2+\frac{1}{\gamma-1}\rho^\gamma)(t,r)r^2\,\mathrm{d} r\\
&\quad+8\pi^2\delta\int_0^\infty\int_0^\infty K(l,r)\rho(t,l)\rho(t,r)l^2r^2 \,\mathrm{d} l \mathrm{d} r,
\end{align*}
where the kernel \(K(l,r)\) is defined by
\begin{eqnarray*}
	K(l,r)=	
	\left\{ \begin{array}{ll}
		1/l \quad \text{if}\ r\leq l,\\
		1/r \quad \text{if}\ r\geq l.
	\end{array}
	\right.
\end{eqnarray*}
We should emphasize that the mass and physical energy are only used in the formation of singularities of Euler-Poisson equations with electrostatic forces (\(\delta=1\)).

Apart from the mass and physical energy, we further introduce  the weighted total radial velocity
\begin{align*}
F(t)=-\int_0^\infty v(t,r)\exp(-r)\, \mathrm{d} r,
\end{align*}
which will play a crucial role in studying  formation of singularities afterwards. It is clear that the functional \(F(t)\) is well defined as long as
the solution \((\rho,u)\) stays in \(\mathcal{X}_T\) via Sobolev embedding.

Let \((\rho,u)\in \mathcal{X}_T\)  be the  solution to \eqref{1} or \eqref{3} with  the spherically symmetric initial data \eqref{10.5}.

We first state singularity formation  of the Euler equations (\(\delta=0\)).
\begin{Theorem}\label{theorem:Euler}  If the initial data \((\rho_0,u_0)\) satisfy
	\begin{align*}
	\rho_0(0)=0,\ F(0)>0.
	\end{align*}
	Then,  the lifespan \(T^*\) of  the classical solution \((\rho,u)\) to \eqref{1} with \eqref{10.5}  is bounded above by
	\begin{align}\label{14}
	T^*\leq\frac{2}{F(0)}.
	\end{align}
\end{Theorem}
\textbf{Proof.}
We will prove that the lifespan \(T^*\) is finite by contradiction. Assume \(T^*=\infty\).
It follows from the momentum equation of \eqref{12} that
\begin{equation}\label{15}
\begin{aligned}
\frac{\mathrm{d}}{\mathrm{d} t}F(t)
&=\int_0^\infty\big[v\partial_rv+\frac{\gamma}{\gamma-1}\partial_r(\rho^{\gamma-1})\big]\exp(-r)\,\mathrm{d} r.
\end{aligned}
\end{equation}
Integrating by parts in view of \eqref{12.5}, one has
\begin{align}\label{16}
\int_0^\infty v\partial_rv\exp(-r)\,\mathrm{d} r
=\frac{1}{2}\int_0^\infty v^2\exp(-r)\,\mathrm{d} r
\end{align}
and
\begin{align}\label{17}
\int_0^\infty\frac{\gamma}{\gamma-1}\partial_r(\rho^{\gamma-1})\exp(-r)\,\mathrm{d} r
=\frac{\gamma}{\gamma-1}\int_0^\infty \rho^{\gamma-1}\exp(-r)\,\mathrm{d} r.
\end{align}
Inserting \eqref{16} and \eqref{17} into \eqref{15} yields   	
\begin{equation}\label{18}
\begin{aligned}
\frac{\mathrm{d}}{\mathrm{d} t}F(t)
&= \frac{1}{2}\int_0^\infty v^2\exp(-r)\,\mathrm{d} r+\frac{\gamma}{\gamma-1}\int_0^\infty \rho^{\gamma-1}\exp(-r)\,\mathrm{d} r\\
&\geq \frac{1}{2}\int_0^\infty v^2\exp(-r)\,\mathrm{d} r=\colon Q(t),
\end{aligned}
\end{equation}	
in which we have used the fact that \(\rho\) is non-negative in the last inequality.

To control \(F(t)\) in terms of \(Q(t)\),  we employ Cauchy-Schwartz inequality  to deduce
\begin{equation}\label{19}
\begin{aligned}
F^2(t)
\leq\int_0^\infty v^2\exp(-r)\,\mathrm{d} r\int_0^\infty\exp(-r)\,\mathrm{d} r=2Q(t).
\end{aligned}
\end{equation}

We conclude from \eqref{18} and \eqref{19} that
\begin{align*}
\frac{\mathrm{d}}{\mathrm{d} t}F(t)\geq \frac{1}{2}F^2(t).
\end{align*}
It follows that
\begin{align*}
F(t)\geq \frac{F(0)}{1-\frac{F(0)}{2}t}.
\end{align*}
This implies that \(F(t)\) will develop singularity no later than the time \(\frac{2}{F(0)}\). This contradiction ends the proof.
\qed

\begin{Remark} The Euler equations in general dimension may be written in spherical coordinates as   
\begin{eqnarray*}
\left\{ \begin{array}{ll}
\partial_t\rho+\frac{1}{r^{n-1}}\partial_r(r^{n-1}\rho
v)=0,\\
\rho(\partial_tv+v\partial_rv)+\partial_rp=0.
\end{array}
\right.
\end{eqnarray*}
The proof of Theorem \ref{theorem:Euler} can be applied  uniformly here and the same result follows. Our proof also simplifies that of \cite{Lei2013singularities}.    

\end{Remark}

The following is the formation singularity result for the Euler-Poisson equations in the presence of gravitational forces (\(\delta=-1\)).
\begin{Theorem}\label{theorem:EPE}   If  the initial data \((\rho_0,u_0)\) satisfy
	\begin{align*}
	\rho_0(0)=0,\ F(0)>0.
	\end{align*}
	Then,  the lifespan \(T^*\) of  \((\rho,u)\) to \eqref{3} (\(\delta=-1\)) with  \eqref{10.5} is bounded above by
	\begin{align*}
	T^*\leq\frac{2}{F(0)}.
	\end{align*}
\end{Theorem}
\textbf{Proof.}
One calculates from the momentum equation of \eqref{12} that
\begin{equation*}
\begin{aligned}
\frac{\mathrm{d}}{\mathrm{d} t}F(t)
&=\int_0^\infty\bigg[v\partial_rv+\frac{\gamma}{\gamma-1}\partial_r(\rho^{\gamma-1})+\frac{4\pi }{r^2}\int_0^r\rho(t,l) l^2\,\mathrm{d}l\bigg]\exp(-r)\,\mathrm{d} r\\
&\geq\int_0^\infty\bigg[v\partial_rv+\frac{\gamma}{\gamma-1}\partial_r(\rho^{\gamma-1})\bigg]\exp(-r)\,\mathrm{d} r.
\end{aligned}
\end{equation*}
The rest is same as the proof of Theorem \ref{theorem:Euler}.    	
\qed

In contrast to the gravitational forces case, without good sign to use,  the blowup formation in the case of the electrostatic forces (\(\delta=1\)) is slightly more complicated.
\begin{Theorem}\label{theorem:EPG} Let \(\sigma\in (0,\frac{1}{2})\) and \(\frac{3-\sigma}{2-\sigma}<\gamma\leq \frac{5}{3}\).
	If   \((\rho_0,u_0)\) satisfy
	\begin{align*}
	\rho_0(0)=0
	\end{align*}	
	and
	\begin{align}\label{20}
	\frac{1}{4}F^2(0)\geq
	c_1\bigg(M(0)+(\gamma-1)E_1(0)\bigg)+c_2,
	\end{align}
	where
	\begin{align*}
	c_1=\frac{2-\sigma}{3-\sigma}\big(\frac{3-\sigma}{3-3\sigma+\sigma^2}\big)^{\frac{3-\sigma}{2-\sigma}},\ c_2=\frac{4\pi e^{-1}(1+\sigma)}{\sigma(3-\sigma)}.
	\end{align*}	
	Then  the lifespan \(T^*\) of  \((\rho,u)\) to \eqref{3} (\(\delta=1\)) with  \eqref{10.5} is bounded above by
	\begin{align*}
	T^*\leq\frac{4}{F(0)}.
	\end{align*}
\end{Theorem}
\textbf{Proof.}
It follows from the momentum equation of \eqref{12} that
\begin{equation}\label{21}
\begin{aligned}
\frac{\mathrm{d}}{\mathrm{d} t}F(t)
&=\int_0^\infty\bigg[v\partial_rv+\frac{\gamma}{\gamma-1}\partial_r(\rho^{\gamma-1})-\frac{4\pi }{r^2}\int_0^r\rho(t,l) l^2\,\mathrm{d}l\bigg]\exp(-r)\,\mathrm{d} r.
\end{aligned}
\end{equation}
Denote
\begin{align*}
R(t)&=4\pi\int_0^\infty\frac{1}{r^2}\int_0^r\rho(t,l) l^2\,\mathrm{d}l\exp(-r)\,\mathrm{d} r.
\end{align*}
There is too much singularity near \(r=0\) in the integrand of \(R(t)\), which can not be handled by  weighted Hardy inequality directly.
To get round  this difficulty, by compensating a  power of \(r^\sigma\), via Young inequality and weighted Hardy's inequality from Appendix, one has
\begin{equation}\label{22}
\begin{aligned}
R(t)&=4\pi\int_0^\infty\frac{1}{r^2}\int_0^r\rho(t,l) l^{2-\sigma}l^\sigma\,\mathrm{d}l\exp(-r)\,\mathrm{d} r\\
&\leq 4\pi\int_0^\infty\frac{1}{r^{2-\sigma}}\int_0^r\rho(t,l) l^{2-\sigma}\,\mathrm{d}l\exp(-r)\,\mathrm{d} r\\
&\leq 4\pi\int_0^\infty\frac{1}{r^{1-\sigma}}\bigg[{\frac{2-\sigma}{3-\sigma}}\left(\frac{1}{r}\int_0^r\rho(t,l) l^{2-\sigma}\,\mathrm{d}l\right)^{\frac{3-\sigma}{2-\sigma}}+{\frac{1}{3-\sigma}}\bigg]\exp(-r)\,\mathrm{d} r\\
&\leq \frac{4\pi(2-\sigma)}{3-\sigma}\int_0^\infty\frac{1}{r^{1-\sigma+\frac{3-\sigma}{2-\sigma}}}\left(\int_0^r\rho(t,l) l^{2-\sigma}\,\mathrm{d}l\right)^{\frac{3-\sigma}{2-\sigma}}\,\mathrm{d} r+
\frac{4\pi}{3-\sigma}\int_0^\infty\frac{\exp(-r)}{r^{1-\sigma}}\,\mathrm{d} r\\
&\leq \frac{4\pi(2-\sigma)}{3-\sigma}\big(\frac{3-\sigma}{3-3\sigma+\sigma^2}\big)^{\frac{3-\sigma}{2-\sigma}}\int_0^\infty
\rho^{\frac{3-\sigma}{2-\sigma}}(t,r)r^2\,\mathrm{d} r+
\frac{4\pi e^{-1}(1+\sigma)}{\sigma(3-\sigma)},
\end{aligned}
\end{equation}
and the first term in \eqref{22} may be continuously estimated as
\begin{equation}\label{23}
\begin{aligned}
&\int_0^\infty\rho^{\frac{3-\sigma}{2-\sigma}}(t,r)r^2\,\mathrm{d} r\\
&=\int_{(0,\infty)\cap\{\rho(t,r)\leq 1\}}
\rho^{\frac{3-\sigma}{2-\sigma}}(t,r)r^2\,\mathrm{d} r+\int_{(0,\infty)\cap\{\rho(t,r)>1\}}
\rho^{\frac{3-\sigma}{2-\sigma}}(t,r)r^2\,\mathrm{d} r\\
&\leq\int_{(0,\infty)\cap\{\rho(t,r)\leq 1\}}
\rho(t,r)r^2\,\mathrm{d} r+\int_{(0,\infty)\cap\{\rho(t,r)>1\}}
\rho^\gamma(t,r)r^2\,\mathrm{d} r\\
&\leq\int_0^\infty
\rho(t,r)r^2\,\mathrm{d} r+\int_0^\infty
\rho^\gamma(t,r)r^2\,\mathrm{d} r\\
&\leq \frac{1}{4\pi}\big[M(0)+(\gamma-1)E_1(0)\big],
\end{aligned}
\end{equation}
where we have used \(\frac{3-\sigma}{2-\sigma}<\gamma\leq \frac{5}{3}\).  Substituting \eqref{23} into \eqref{22} gives
\begin{align}\label{24}
R(t)\leq c_1\bigg(M(0)+(\gamma-1)E_1(0)\bigg)+c_2,
\end{align}
in which \(c_1=\frac{2-\sigma}{3-\sigma}\big(\frac{3-\sigma}{3-3\sigma+\sigma^2}\big)^{\frac{3-\sigma}{2-\sigma}}\) and \(c_2=\frac{4\pi e^{-1}(1+\sigma)}{\sigma(3-\sigma)}\).
It follows from \eqref{21} and \eqref{24} that
\begin{align}\label{25}
\frac{\mathrm{d}}{\mathrm{d} t}F(t)\geq Q(t)-
c_1\bigg(M(0)+(\gamma-1)E_1(0)\bigg)-c_2.
\end{align}

We conclude from \eqref{19} and \eqref{25} that
\begin{align*}
\frac{\mathrm{d}}{\mathrm{d} t}F(t)\geq \frac{1}{2}F^2(t)-
c_1\bigg(M(0)+(\gamma-1)E_1(0)\bigg)-c_2.
\end{align*}
In view of \eqref{20}, by a bootstrap argument one shows that \(F(t)\) is increasing and thus satisfies
\begin{align*}
F(t)\geq \frac{F(0)}{1-\frac{F(0)}{4}t},
\end{align*}
which shows that \(F(t)\) will breakdown before the time \(\frac{4}{F(0)}\). This is a contradiction.

\qed

\section{Sharp Bound and Example}
Let \((\rho,u)\in \mathcal{X}_T\)  be  solution to \eqref{1} or \eqref{3} with  the spherically symmetric initial data \eqref{10.5}. Based on the local well-posedness in Section 3 and formation of singularities in Section 4 we can obtain a sharp bound on the lifespan of the solutions to the Euler equations and  Euler-Poisson equations with gravitational forces by choosing a special scale relation between the initial density and initial velocity. We also give an example to justify the conditions of the formation of singularity of the Euler-Poisson equations with electrostatic forces.

We have the following  sharp bound on lifespan of the solutions to the Euler equations and  Euler-Poisson equations with gravitational forces.
\begin{Theorem} Assume that the initial density \(\rho_0\) and velocity \(u_0\) satisfy the scaling relation  \((\rho_0,u_0)=(\varepsilon^\frac{2}{\gamma-1}\varrho, \varepsilon w)\). If the initial data \((\rho_0,u_0)\) also satisfy
	\begin{align}\label{25.5}
	\rho_0(0)=0,\ F(0)>0.
	\end{align}
	Then, the lifespan \(T^*\) of  \((\rho,u)\) to \eqref{1} with \eqref{10.5} satisfies
	\begin{align}\label{26}
	T^*\eqsim\frac{1}{\varepsilon}.
	\end{align}
\end{Theorem}
\textbf{Proof.} A direct calculation gives
\begin{align*}
\mathcal{E}(0)\eqsim\varepsilon^2, \ F(0)\eqsim\varepsilon.
\end{align*}
The conclusion \eqref{26} follows from  \eqref{5} and  \eqref{14}.
\qed
\begin{Theorem} Assume that the initial density \(\rho_0\) and velocity \(u_0\) satisfy the scale relation  \((\rho_0,u_0)=(\varepsilon^\frac{2}{\gamma-1}\varrho, \varepsilon w)\). If the initial data \((\rho_0,u_0)\) also satisfy
	\begin{align*}
	\rho_0(0)=0,\ F(0)>0.
	\end{align*}
	Then, the lifespan \(T^*\) of  \((\rho,u)\) to \eqref{3} (\(\delta=-1\)) with  \eqref{10.5} satisfies
	\begin{align*}
	T^*\eqsim\frac{1}{\varepsilon}.
	\end{align*}
\end{Theorem}

\begin{Remark} Let
	\begin{align}\label{27}
	\varrho(x)=|x|^2\exp(-|x|^2), \  w(x)=-x|x|^2\exp(-|x|^2).
	\end{align}
	It is clear that \((\rho_0,u_0)\) belong to \(H^3(\mathbb{R}^3)\) and satisfy \eqref{25.5}. 	
\end{Remark}

We finally give an example to justify \eqref{20}. Taking \((\varrho, w)\) of \eqref{27} and setting \((\rho_0,u_0)=(\varepsilon\varrho, \varepsilon^{-1} w)\), one calculates that
\begin{align*}
M(0)\eqsim\varepsilon,\quad   F(0)\eqsim\varepsilon^{-1}
\end{align*}
and
\begin{align*}
E_1(0)\eqsim\varepsilon^{-1}+\varepsilon^{\gamma}-\varepsilon^2\eqsim\varepsilon^{-1}.
\end{align*}
Consequently, for given \(\sigma\in (0,\frac{1}{2})\), let \(\varepsilon\) go to zero, one has 	
\begin{align*}
\frac{F^2(0)}{c_1\bigg(M(0)+(\gamma-1)E_1(0)\bigg)+c_2}\eqsim\varepsilon^{-1}\gg 1.
\end{align*} 	
This implies that \eqref{20} holds.

\section{Appendix}

We state here the Hardy-Littlewood-Sobolev inequality and weighted Hardy inequality for references convenience,  whose proofs can be found  in \cite{stein2016singular} and \cite{kufner2007hardy}, respectively.
\begin{Lemma} Given \(0<\lambda<d\) and \(1<q_1<q_2<\infty\) with \(\frac{1}{q_2}=\frac{1}{q_1}-\frac{\lambda}{d}\). Let \(I_{\lambda}\)  be the Riesz potential of order \(\lambda\) on \(\mathbb{R}^d\) which is defined by
	\begin{align*}
	I_{\lambda}f(x)=C(d,\lambda)\int_{\mathbb{R}^d}\frac{f(y)}{|x-y|^{d-\lambda}}\,\mathrm{d} y.
	\end{align*} 	
	Then
	\begin{align*}
	\|I_{\lambda}f\|_{L^{q_2}(\mathbb{R}^d)}\leq C(d,\lambda,q_1)\|f\|_{L^{q_1}(\mathbb{R}^d)}.
	\end{align*}	
\end{Lemma}
\begin{Lemma} Given  \(1<\mu,q<\infty\).
	Then
	\begin{align*}
	\int_0^\infty\frac{1}{r^\mu}\left(\int_0^rf(l)\,\mathrm{d}l\right)^q\,\mathrm{d} r\leq C(\mu,q) \int_0^\infty f^q(r)r^{q-\mu}\,\mathrm{d} r,
	\end{align*}
	where the best constant \(C(\mu,q)=\left(\frac{q}{\mu-1}\right)^q\). In particular, it is the classical Hardy's inequality as \(\mu=q\).	
\end{Lemma}

\vskip 5mm

\noindent \textrm{\textbf{Acknowledgements:}} The research of Li was supported partially by the National Natural Science Foundation of China (Nos. 11231006, 11225102, 11461161007 and 11671384), and the Importation and Development of High Caliber Talents Project of Beijing Municipal Institutions (No. CIT\&TCD20140323).
The research of Wang was supported by grant nos. 231668 and 250070 from the Research Council of Norway.

\noindent \textrm{\textbf{Conflict of Interest:}} The authors declare that they have no conflict of interest.


\textsc{School of Mathematics, Capital Normal University, Beijing 100048, P. R. China}

 E-mail address: hailiang.li.math@gmail.com

\textsc{Department of Mathematical Sciences, Norwegian University of Science and Technology, Trondheim 7491, Norway}

E-mail address: yuexun.wang@ntnu.no

\end{document}